\documentclass[11pt]{amsart}
\usepackage{a4wide}
\usepackage{enumerate}
\usepackage{amssymb}
\usepackage{amsthm}
\usepackage{exscale}
\newcommand{\RR}{\mathbb{R}}
\newcommand{\CC}{\mathbb{C}}

\newcommand{\NN}{\mathbb{N}}
\newcommand{\ZZ}{\mathbb{Z}}

\newcommand{\EE}{\, \mathbb{E} \,}

\newcommand{\Tr}{{\mathop{\mathrm{Tr} \,}}}
\newcommand{\Var}{BV}
\newcommand{\per}{\mathrm{per}}
\renewcommand{\L}{\Lambda}
\newcommand{\tL}{Q}

\newcommand{\supp}{{\mathop{\mathrm{supp\,}}}}
\newcommand{\be}{\begin{equation}}
\newcommand{\ee}{\end{equation}}

\newtheorem{thm}{Theorem}
\newtheorem{lem}[thm]{Lemma}
\newtheorem{prp}[thm]{Proposition}

\theoremstyle{definition}

\theoremstyle{remark}
\newtheorem{rem}[thm]{Remark}

\newcommand{\hm}[1]{\leavevmode{\marginpar{\tiny%
$\hbox to 0mm{\hspace*{-0.5mm}$\leftarrow$\hss}%
\vcenter{\vrule depth 0.1mm height 0.1mm width \the\marginparwidth}%
\hbox to
0mm{\hss$\rightarrow$\hspace*{-0.5mm}}$\\\relax\raggedright #1}}}

\begin{document}

\title[Wegner estimates for sign-changing single site potentials]
{Wegner estimates for \\ sign-changing single site potentials}

\author[I.~Veseli\'c]{Ivan Veseli\'c}
\address{Emmy-Noether-Programme of the Deutsche Forschungsgemeinschaft \& Fakult\"at f\"ur Mathematik,\, 09107\, TU-Chemnitz, Germany  }
\urladdr{www.tu-chemnitz.de/mathematik/enp/}

\thanks{
{\today, \jobname.tex}}


\begin{abstract}
We study Anderson and alloy type random Schr\"odinger operators on $\ell^2(\ZZ^d)$ and $L^2(\RR^d)$. 
Wegner estimates are bounds on the average number of eigenvalues in an energy interval of finite box restrictions of 
these types of operators.
For a certain class of models we prove a Wegner estimate which is linear in the volume of the box
and the length of the considered energy interval. The single site potential of the Anderson/alloy type model does not need to have fixed sign,
but it needs be of a generalised step function form.
The result implies the Lipschitz continuity of the integrated density of states.
\end{abstract}

\maketitle

\section{Model and results}

We study spectral properties of Schr\"odinger operators which are given as the sum $H=-\Delta +V$
of the negative Laplacian $\Delta $ and a multiplication operator $V$. 
The operators can be 
considered in $d$-dimensional Euclidean space $\RR^d$ or on the lattice $\ZZ^d$. 
To be able to treat both cases simultaneously let us use the symbol $X^d$ for either 
$\RR^d$ or $\ZZ^d$. 
On the continuum the Laplace operator is the sum of second derivatives	 $\sum_{i=1}^d \frac{\partial^2}{\partial x_i^2}$
and $V$ is a bounded function $\RR^d \to \RR$. Thus $H$ is selfadjoint on the usual Sobolev space $W^{2,2}(\RR^d)$.
In the discrete case the Laplacian is given by the rule $\Delta \phi(k)= \sum_{i=1}^d \phi(k + e_i)+ \phi(k - e_i)$,
where $\phi$ is a sequence in $\ell^2(\ZZ^d)$ and $(e_1, \dots, e_d )$ is an orthonormal basis which defines the lattice $\ZZ^d$ as a subset of $\RR^d$.
The potential is given by a bounded function $V \colon \ZZ^d\to \RR$, and thus $H$ is a bounded selfadjoint operator.

The operators we are considering are random. More precisely, the potential $V=V_\per+V_\omega$ decomposes 
into a part $V_\per$ which is translation invariant with respect to some sub-lattice $n\ZZ^d$, $n \in\NN$, i.e.{}
$V_\per(x+k)=V_\per(x)$ for all $ x\in \RR^d$ and all $k \in n\ZZ^d$, and a part $V_\omega$ which is random.
The random part of the potential is a stochastic field 
$V_\omega(x) := \sum_{k \in\ZZ^d} \omega_k u(x-k), x \in X^d$, of \textit{alloy} or \textit{Anderson type}. Here $u \colon X^d \to  \RR$ 
is a bounded, compactly supported function, which we call single site potential. The 
coupling constants $\omega_k, k \in\ZZ^d$ form an independent, identically distributed sequence of real  random variables. 
We assume that the random variables are bounded and distributed according to a density $f$ of bounded variation.
In the discrete case the random operator $H_\omega= -\Delta + V_\per + V_\omega$ is called \textit{Anderson model}, 
and in the continuum case $H_\omega$ is called \textit{alloy type model}.

There is a well defined spectral distribution function
$ N \colon \RR \to \RR$ of the family $(H_\omega)_\omega$ which 
is closely related to eigenvalue counting functions on finite cubes.
To explain this precisely, we need some more notation.
Denote by $\chi$ the characteristic function of  the set
$[-1/2,1/2]^d \, \cap \, X^d$. Thus in the continuum case this set is a unit cube,
and in the discrete case it is a single point. 
Also, for $ k \in \ZZ^d$, let $\chi_k (x): = \chi(x-k)$ be the translate of $\chi$. 
The cube $[-l -\frac{1}{2}, l + \frac{1}{2}]^d \, \cap \, X^d$  will be abbreviated by $\Lambda_l$,
its intersection with $\ZZ^d$ by $\tL_l$,
the restriction of $H_\omega$ to $\Lambda_l$  with selfadjoint boundary conditions 
(e.g.~Dirichlet, Neumann or periodic ones) by $H_\omega^l$, 
the spectral projection associated to $H_\omega$ (respectively to $H_\omega^l$) 
and an interval $I$ by $P_\omega(I)$ (respectively by $P_\omega^l(I)$),
and the number of eigenvalues of $H_\omega^l$
in $]-\infty,E]$  by $N_\omega^l(E):=\mathop{\mathrm{Tr}} [P_\omega^{\Lambda_l}(\,]-\infty,E])]$. 
With  this notation we can define the \textit{integrated density of states},
which is the spectral distribution of the family $(H_\omega)_\omega$,
by  
\[
 N(E) := \EE \{\mathop{\mathrm{Tr}} [\chi P_\omega(]-\infty, E])]\}. 
\]
Here $\chi$ is understood as a multiplication operator.
The function $N$ has the following self-averaging property: for all $E$ where $N$ is continuous (that's a set with countable complement)
the relation $\lim_{l\to\infty}  (2l+1)^{-d} \, N_\omega^l(E) =N(E)$
holds almost surely. This implies that if $E_1$ and $E_2$ are two continuity points of $N$, 
we have 
\begin{equation}
\label{e:1}
\lim_{l\to\infty}  (2l+1)^{-d}\,  \EE\{ N_\omega^l(E_2)-N_\omega^l(E_1)\} = 
N(E_2)-N(E_1).
\end{equation}
Thus if one is able to show that there is a function $C \colon \RR \to \RR$
and an exponent $ \beta \in \, ] 0,1]$, 
such that for all $ E_1, E_2 \le E$ and for all $ l \in \NN$ the so-called Wegner bound
(named after the paper \cite{Wegner-81})
\begin{equation}
\EE\{N_\omega^l(E_2)-N_\omega^l(E_1)\} \le C_W(E) \,  (2l+1)^{d} \, | E_2-E_1|^\beta 
\end{equation}
holds, it follows that the integrated density of states  is (locally uniformly) H\"older-continuous with exponent $\beta$.
Note that this shows a posteriori, that there are no points of discontinuity of $N$ and thus
the convergence in \eqref{e:1} hold actually for all $E_1,E_2\in\RR$. 
This is only one of the reasons why one is interested in bounds on the averaged quantity $ \EE \{\mathop{\mathrm{Tr}} [P_\omega^{\Lambda_l}(]E_1,E_2])]\}$. 
It plays also a crucial role in arguments leading to the proof of localisation, i.e.~the phenomenon that 
there is a subset $I_{loc}\subset \RR$  such that $I_{loc}\, \cap \, \sigma_{pp}(H_\omega)=I_{loc}$
and  $I_{loc}\, \cap \, (\sigma_{ac}(H_\omega)\, \cup \,\sigma_{sc}(H_\omega))=\emptyset$ almost surely. In fact,
usually localisation goes along with quite explicit bounds on the decay of eigenfunctions
and on the non-spreading of electron wavepackets (see for instance the monograph \cite{Stollmann-01} 
or the characterisation established in \cite{GerminetK-04}). For recent surveys on the integrated density of states
see \cite{KirschM-07,Veselic-07b}.\\

Now we specialise to a specific class of single site potentials. 
The important point is that the resulting Anderson/alloy type model allows for
random potentials where the potential values at different points in space are
negatively correlated: Let $ \kappa >0$, $v \colon X^d \to \RR$ a function satisfying $v \ge \kappa \chi$,
and $\alpha \colon \ZZ^d\to\RR$ a function with compact support 
such that its Fourier transform $\hat\alpha \colon [0, 2\pi \, [^d \to \CC$, 
$\hat\alpha (\theta) :=\sum_{k \in \ZZ^d} \alpha_k \mathrm e^{-\mathrm i k \cdot \theta}$
does not vanish on $[0, 2\pi \, [^d $.
 Then we call 
\begin{equation}
\label{e:3}
u(x):= \sum_{k \in \ZZ^d} \alpha_k v (x-k)
\end{equation}
a single site potential of \textit{generalised step function form}. Note that the 
sum contains only finitely many non-vanishing terms. Due to the fact that the coefficients
$\alpha_k, k \in \ZZ^d,$ may change sign, the random potential $V_\omega$ 
can have negative correlations between values at different sites. 
Now we are in the position to formulate 
our main result:
\smallskip

\begin{thm}
\label{t:main}
Let $H_\omega= -\Delta+ V_\per + V_\omega$ be an Anderson model on $\ell^2(\ZZ^d)$
or an alloy type model on $L^2(\RR^d)$ with a single site potential $u$ 
of generalised step function form. Assume that the density $f$
of the coupling constants has compact support and bounded variation.
Then there is a continuous function $C_W \colon \RR \to \RR$
such that for all $ E_1, E_2 \le E$ and for all $ l \in \NN$ the Wegner bound
\begin{equation}
\label{e:4}
\EE\{N_\omega^l(E_2)-N_\omega^l(E_1)\} \le C_W(E) \,  (2l+1)^{d} \, | E_2-E_1| 
\end{equation}
holds. 
 \end{thm}
A possible choice of the function $C_W$ is 
$C_W(E)= \frac{1}{\kappa}C(E,V) \, \|f\|_{\Var} \, \|B\|_1$.
Here $\|f\|_{\Var}$ denotes the total variation norm of $f$,
$C(E,V) := e^{E+ V_\infty} \, \sum_{n \in \ZZ^d, n_j \ge 0} \exp\big (\frac{\pi^2}{2} \sum_{j=1}^d n_j^2 \big )$,
$V_\infty:= \sup_{x \in X, \omega\in\Omega} |V_\per(x)+ V_\omega(x)|$, and $\|B\|_1$ is the column sum norm of
the inverse of the multi-dimensional Laurent matrix $\{\alpha_{j-k}\}_{j,k\in\ZZ^d}$.
In the case of the discrete Anderson model one can choose $C_W(E)= \frac{1}{\kappa} \, \|f\|_{\Var} \, \|B\|_1$.
\\

To apply the theorem one has to know that $\hat\alpha$ vanishes nowhere on the torus.
Let us give two instances where this condition holds.
The first case is when there is an index $j \in \ZZ^d$ such that 
$|\alpha_j| > \sum_{k \in \ZZ^d, k \neq j} |\alpha_k|$. Secondly, in the case that $d=1$ 
and that the diameter $N$ of the support of $\alpha$ is kept fixed, 
the property holds for a dense, open subset of $\alpha \in\RR^{N+1}$.

As mentioned above, estimate \eqref{e:4} implies that the integrated density of states  $N\colon \RR \to\RR$ 
is locally uniformly Lipschitz-continuous. This in turn implies that the 
derivative $ n(E):=\frac{dN(E)}{dE}$ exists almost everywhere on $\RR$ and is 
locally uniformly bounded by $ n(E) \le C_W(E)$. The function $n$ is called \textit{density of states}.
\\

For Anderson/alloy type models where the single site potential has fixed sign, 
Wegner estimates  are well understood by now, see e.g.{} \cite{KirschM-07,CombesHK-07,Veselic-07b}.
Let us discuss earlier theorems in the literature which establish Wegner estimates for 
single site potentials that change sign.

Theorem \ref{t:main} recovers the main result of \cite{Veselic-02a}
where the same statement was proven under two additional conditions:
It was assumed that there is an index $j \in \ZZ^d$ such that 
$|\alpha_j| > \sum_{k \in \ZZ^d, k \neq j} |\alpha_k|$ and that 
the density $f$ belongs to the Sobolev space $W^{1,1}_c(\RR)$.
Exactly the same statement as in Theorem \ref{t:main} above, but only for 
dimensions $d=1$ and $d=2$ was proven in \cite{KostrykinV-06} in a joint paper with V.~Kostrykin.
There is another method to prove Wegner estimates for 
single site potentials that are allowed to change sign which is based 
on certain vector fields in the parameter space underlying the alloy type model.
It was introduced in \cite{Klopp-95a} by F.~Klopp 
and improved by P.~Hislop and F.~Klopp  in \cite{HislopK-02}.
Its advantage is that it applies to arbitrary continuous, compactly
supported single site potentials (which are not identically equal to zero).
The regularity requirement on the density $f$ is slightly more restrictive
than in the Theorem \ref{t:main}.
However, this method applies only to certain energy
intervals $[E_1,E_2]$: sufficiently low energies are allowed, 
but arbitrary high energies are not allowed.
The papers  \cite{Klopp-95a,HislopK-02,KostrykinV-06} 
contain various other results, which we do not state here, 
because they cannot be directly compared with our theorem above.
Additional aspects of Wegner estimates for sign-nondefinite single site potentials are 
discussed e.g. in \cite{CombesHKN-02}, Section 5.5. of \cite{Veselic-07b}, and \cite{CombesHK-07}.
\smallskip

Let us briefly discuss the relevance of the condition that the Fourier
transform $\hat \alpha$ does not vanish on $[0,2\pi\,[^d$. It ensures that 
the multi-dimensional Laurent matrix $A$ with coefficients 
$\alpha_{j-k}, j,k\in\ZZ^d$, when considered as an operator 
from $\ell^p(\ZZ^d)$  to $\ell^p(\ZZ^d)$ has a bounded inverse $B$.
However, in the proof of the Wegner estimate above we 
encounter not the infinite matrix $A$, but rather finite
size matrices $A_\Lambda$ which need to have bounded inverses $B_\Lambda$
with norms uniformly bounded in $\Lambda=\Lambda_l, l \in \NN$.
The relevant norm is the column sum norm, 
corresponding to the operator norm on $\ell^1(\Lambda)$.
If  $A_\Lambda$ is chosen to be a finite section multi-dimensional
Toeplitz operator this leads to nontrivial open questions concerning 
the invertibility of truncated Toeplitz matrices, see for instance 
\cite{BoettcherS-99}. This is the reason 
why the results of \cite{KostrykinV-06} are restricted to dimension one and two,
cf.~also \cite{KozakS-80b}. However, it turns out that one has a certain freedom 
in the choice of the finite volume matrices $A_\Lambda$. In particular, one can 
choose them to be finite multi-dimensional circulant matrices (rather than finite Toeplitz 
matrices), which have much better invertibility properties and can be used to 
complete the proof of Theorem \ref{t:main}.

\section{Proof of Theorem \ref{t:main}}

\begin{rem}[Cubes]
Theorem \ref{t:main} concerns Hamiltonians restricted to a cube $\Lambda_l\subset X^d$ of side length $2l+1$.
However, in the proof we will have to deal with several modifications of this cube.
First, we need to consider cubes which are subsets of $\ZZ^d$ and not of $X^d$ (which may be either $\RR^d$ or  $\ZZ^d$).
Second, modified cubes will be of larger side length than the original cube $\Lambda_l$. 
At this point we list the various cube sizes which will appear at various stages of the proof.

To avoid confusion, cubes in $X^d$ will be denoted by the symbol $\Lambda$ 
while the letter $\tL$ will be reserved for cubes in $\ZZ^d$. 
In the following we assume that $l \in \NN$ is fixed and that $\Lambda= \Lambda_l$.
Since the single site potential $u$ has compact support, there is some $R\in\NN$ such 
that $\supp u $ is contained in $\Lambda_R$. This implies that the set of 
the lattice sites $k$ such that the coupling constants $\omega_k$ influence the 
potential $V_\omega$ inside the cube $\Lambda$ is contained in the set $\tL_{l+R}=\L_{l+R} \cap \ZZ^d$.
Similarly, there is some $r\in\NN$ such that  the support of $v$ is contained in $\Lambda_r$. Consequently, 
the set $\{ k\in\ZZ^d\mid  \supp v_k  \cap \Lambda_l \}$ is contained in $\tL_{l+r}$.
Here we used the abbreviation $v_k(x)=v(x-k)$.
Finally, there is some $D\in \NN$ such that $\supp\alpha\subset Q_D$. The relation between $v,u$, and $\alpha$
implies $R \le r+D$. Let us point point out that all cubes listed depend on the reference cube $\Lambda$ of side length $2l+1$.
\end{rem}

The proof of Theorem \ref{t:main} uses results on spectral averaging for \emph{non-negative} single site potentials
established  in \cite[\S 4]{CombesH-94b} (see also \cite{KotaniS-87}).
These facts are formulated in Proposition \ref{p:CombesHislop}.
The subsequent Propositions \ref{p:main1} and \ref{p:main2} contain the estimates which are
needed to deal with single site potentials \emph{of changing sign}

\begin{prp}[\cite{CombesH-94b}] \label{p:CombesHislop} 
Let $I=[E_1,E_2]$ be an interval. Then
\begin{enumerate}[(a)]
 \item 

\label{TrLoc}
    \[
    \EE \left \{ \Tr P_\omega^l(I) \right \}
    \le C(E_2,V)\, \sum_{j\in \tL_l} \left \|  \EE \{\chi_j P_\omega^l(I) \chi_j  \} \right \|
    \] 
Here $C(E_2,V) := e^{E_2+ V_\infty} \, \sum_{n \in \ZZ^d, n_j \ge 0} \exp\big (\frac{\pi^2}{2} \sum_{j=1}^d n_j^2 \big )$
and $V_\infty:= \sup_{x \in X, \omega\in\Omega} |V_\per(x)+ V_\omega(x)|$.
\item 
Let $ H= -\Delta+W$ be a Schr\"odinger operator with a bounded potential $W$, 
$w$ a function satisfying $w\ge \chi_j$ for some $j\in\ZZ^d$, 
and $ t\mapsto H_t = H + t w$ a one parameter family of operators.
Denote by $H_t^\Lambda$ a selfadjoint restriction of $H_t$ to a cube $\Lambda$  and
by $P_t^\Lambda(I)$ the associated spectral projection on to the interval $I$.
Then, for any $ g \in L_c^\infty(\RR)$ and any  $\phi \in L^2(\Lambda )$ with $\| \phi\|=1$
we have
\[
\int d t \, g(t) \, \langle \phi , \chi_j P_t^\Lambda(I) \chi_j \phi \rangle 
\le 
|I |  \, \|g\|_\infty 
\]
\end{enumerate}
\end{prp}
The first statement allows one to decompose the expectation value of the trace of the spectral projector
into local contributions of unit cubes. The local contributions do not depend on the trace, but rather on the 
norm of certain restricted operators. Spectral averaging is easier to perform of norms than on traces.
Estimate (a) is proven using Dirichlet-Neumann-bracketing and Jensen's inequality in the case of operators
on $L^2(\RR^d)$.  For operators on $\ell^2(\ZZ^2)$ the estimate is trivial and the constant $C(E_2,V)$ can be chosen equal to one.

Statement (b) is a spectral averaging estimate which is based on a contour integral in the complex 
plane and the residue theorem.
\smallskip

The next two propositions contain the main technical results of the paper.
\smallskip

For a cube $\L=\L_l$ let us denote by $A_\Lambda$ and $B_\Lambda$ 
matrices with coefficients in $\tL_{l+R}$. If $A_\Lambda$ is invertible, $B_\Lambda$ 
will denote its inverse. 
The column sum norm $\sup_k \sum_j |B_\L({j,k})| $ of $B_\Lambda$ will be denoted by $\|B_\Lambda\|_1$ . 


\begin{prp} \label{p:main1} 
Let $I=[E_1,E_2]$ be an interval. 
If there exists an invertible matrix $A_\Lambda\colon \ell^1(\tL_{l+R})\to \ell^1(\tL_{l+R})$
such that 
\begin{equation}
 \label{e:Acondition}
 A_\Lambda({j,k}) = \alpha_{j-k}  \quad \text{ for all } j \in \tL_{l+r} \text{ and } k \in \tL_{l+R} 
\end{equation}
then for any $j \in \tL_l$ and $\phi \in L^2(\Lambda_l )$ with $\| \phi\|=1$ 
    \[
    \EE \left \{ \langle \phi , \chi_j P_\omega^l(I) \chi_j \phi \rangle  \right\}
    \le  |I| \, \|f\|_{\Var} \|B_\Lambda\|_1
    \] 
Here $\|f\|_{\Var} $ denotes the total variation norm of the density $f$ of the random variables $\omega_k$.
\end{prp}
Let us point out that condition \eqref{e:Acondition} fixes the coefficients of $A_\L$ only in 
a multi-dimensional \textit{rectangle}. The coefficients outside the rectangle are arbitrary, up to the invertibility condition.
\smallskip

Denote by  $A\colon \ell^1(\ZZ^d)\to \ell^1(\ZZ^d)$
 the linear operator whose coefficients in the canonical orthonormal basis  are 
$A({j,k})= \alpha_{j-k}$ for $j,k \in \ZZ^d$. Since the function $\alpha$ has compact support, 
the operator $A$ is bounded. Moreover, $\hat \alpha$ is an element of the Wiener algebra of functions 
on the $d$-dimensional torus with absolutely convergent Fourier series.
The so-called \mbox{`$1/f$ Theorem'} of Wiener states that if $\hat \alpha$ vanishes nowhere on the torus,
then the inverse $1/{\hat \alpha}$ has an absolutely convergent Fourier series as well.
Wiener's original result concerns the case $d=1$, but Gelfand's proof of the theorem (cf.~e.g.{} \cite{Arveson-02,Katznelson-04})
extends directly to arbitrary $d$.
If we denote by $\beta_n$ the Fourier coefficients of $1/\hat\alpha$, then the Laurent matrix
$B= (\beta_{j-k})_{j,k\in\ZZ^d}$ has finite column sum norm, i.e.~is bounded as an operator 
$\ell^1(\ZZ^d)\to \ell^1(\ZZ^d)$. Furthermore, since  $1/\hat\alpha$ is the inverse of  $\hat\alpha$,
the operator $B$ is the inverse of $A$.

For a given cube $\L=\L_l$ and a site $m \in \ZZ^d$ we consider the associated lattice 
$\Gamma_{l+R} (m) = m + \Gamma_{l+R} $, where $\Gamma_{l+R} =  (2l + 2R +1)\ZZ^d$, 
and the projection $\pi_{l+R} \colon \ZZ^d \to \Lambda_{l+R}$, $\pi_{l+R} (m)=  \Lambda_{l+R} \,\cap\, \Gamma_{l+R} (m)$.
If there is no danger of confusion we will drop the subscript $l+R$ denoting the period.

\begin{prp} \label{p:main2} 
Let $l > R$ and $\Lambda= \Lambda_l$. Define the matrix $A_\Lambda$ by
$A_{\Lambda}({j,k}) = \alpha_{\pi(j-k)}$ for  $j,k\in\tL_{l+r}$.
Then $A_\Lambda$ satisfies condition \eqref{e:Acondition}, is invertible and 
\[
 \|B_\Lambda \|_1  \le \|B \|_1 <\infty
\]
where $B_\Lambda $ is the inverse of $A_\Lambda$ and $B$ the inverse of $A$. 
\end{prp}
Note that $A_\Lambda$ (and thus its inverse $B_\L$, as well) is a multi-dimensional circulant matrix.\\

Now we prove the two propositions before completing the proof of Theorem \ref{t:main} at the end of this section.
The next proof is an adaptation of results in \cite{Veselic-02a,KostrykinV-06}.

\begin{proof}[Proof of Proposition \ref{p:main1}]
Let us first reduce the model to the case $\kappa=1$. 
Obviously we can write the random potential as 
\[
\sum_{k \in\ZZ^d} \omega_k u(x-k)=  \sum_{k \in\ZZ^d} (\kappa \,\omega_k) \,  \big(\frac{1}{\kappa} \, u(x-k) \big).
\]
Now $\frac{1}{\kappa} \, u(x) =\sum_{j \in \ZZ^d} \alpha_j \big(\frac{1}{\kappa} \, v(x-j) \big)$ where 
by assumption $\frac{1}{\kappa} \, v  \ge \chi$.  The distribution of the random variable $\kappa \,\omega_k$
has the density $x \mapsto h(x):=\frac{1}{\kappa} \, f(x/\kappa)$. The total variation norm of $h$ equals $\frac{1}{\kappa} \, \|f\|_{\Var}$.
Thus we can replace $\kappa$ by one, if we keep in mind that the variation norm of the density gets multiplied by $1/\kappa$.

As pointed out earlier $\EE \left \{ \langle \phi , \chi_j P_\omega^l(I) \chi_j \phi \rangle  \right\}$
depends only on a finite number of random variables $\omega_k$. More precisely if $R$ is such that the compact support 
of $u$ is contained in $\Lambda_R$ then only the coupling constants $\omega_k$ with index $k$ in $\tL_{l+R}$ 
influence the scalar product $ \langle \phi , \chi_j P_\omega^l(I) \chi_j \phi \rangle$. 
Thus we can express the expectation value 
$\EE \left \{ \langle \phi , \chi_j P_\omega^l(I) \chi_j \phi \rangle  \right\} $ as a finite dimensional integral
\begin{equation}
\label{e:multi-integral} 
\int_{\RR^L} d\omega_\L  F(\omega_\Lambda) \langle \phi , \chi_j P_\omega^l(I) \chi_j \phi \rangle.
 \end{equation}
Here $F(\omega_\Lambda)=\prod_{k \in {\tL_{l+R}}} f \left ( \omega_k \right )$ 
is the common density of the random variables $\omega_k$ with index in $\tL_{l+R}$
and $L$ is the cardinality of this index set.

In the sequel it will be convenient to have two alternative representations for the 
random potential $V_\omega$. This will be presented next. For any $x \in \RR$ we have 
\begin{multline}
V_\omega(x)
=\sum_{k \in \ZZ^d} \omega_k u (x-k)
=\sum_{k \in \ZZ^d} \omega_k \sum_{j \in\ZZ^d} v(x-j-k)
=\sum_{m \in \ZZ^d} v(x-m)  \sum_{k \in\ZZ^d} \alpha_{m-k} \omega_k. 
\end{multline}
If $x \in \L_l$, the last sum equals 
$\sum_{m \in \tL_{l+r}} v(x-m)  \sum_{k \in\tL_{l+R}} \alpha_{m-k} \omega_k$.
Here $r\in \NN$ is such that the compact support of $v$ is contained in $\L_r$.
Thus, we can conveniently express $V_\omega(x)$ for   $x \in \L_l$ 
using new random variables $\eta_m :=\sum_{k \in\tL_{l+R}} \alpha_{m-k} \omega_k$, $m \in \tL_{l+r}$
as 
\[
V_\omega(x)=\sum_{m \in\tL_{l+r}} \eta_m v(x-m). 
\]
By assumption there exists an invertible matrix $A_\Lambda\colon \ell^1(\tL_{l+R})\to \ell^1(\tL_{l+R})$
such that for all $m \in \tL_{l+r}$
\begin{equation}
\label{e:etaAomega}
\eta_m = \sum_{k \in\tL_{l+R}} A_\Lambda({m,k}) \, \omega_k 
\end{equation}
We define the random variables $\eta_m $ for $ m \in \tL_{l+R} \setminus \tL_{l+r}$
by requiring that the relation \eqref{e:etaAomega} is true for these indices as well.
\bigskip

If we express the random potential in the $\eta$-variables, the `effective' single site potentials $v_j$ 
satisfy $v_j \ge  \chi_j$, and thus the spectral averaging result of Proposition \ref{p:CombesHislop}
(b) applies:
For any $ g \in L_c^\infty(\RR)$, any $j \in \tL_l$, and  $\phi \in L^2(\Lambda_l )$ with $\| \phi\|=1$
we have 
\begin{equation}
 \label{e:CH-for-v}
\int d \eta_j \, g(\eta_j) \, \langle \phi , \chi_j P_{B_\L \eta_\L}^l(I) \chi_j \phi \rangle 
\le 
\vert I \vert  \, \|g\|_\infty.
\end{equation}
However, the $\eta$-random variables are no longer independent. To understand their dependence we have to  analyse the common
density. It can be compactly written in the form $ k(\eta_\L) = | \det B_\L| \, F(B_\Lambda\eta_\Lambda)$
where $F(\omega_\Lambda)$ is the original common density of the $\omega_k, k \in {\tL_{l+R}}$.
Thus \eqref{e:multi-integral} equals
\begin{equation}
\label{e:transformed-integral}
|\det B_\L| \, \int_{\RR^L} d\eta_\L k(\eta_\L)  \langle \phi , \chi_j P_\omega^l(I) \chi_j \phi \rangle
= |\det B_\L| \, \int_{\RR^{L-1}} d\eta_\L ^{\perp j}  \, \int_{\RR} d\eta_j k(\eta_\L)  \langle \phi , \chi_j P_\omega^l(I) \chi_j \phi \rangle.
\end{equation}
Here we denote by $\eta_\L ^{\perp j} $ the sub-collection of random variables indexed by $\tL_{l+R} \setminus \{j\}$.
If we apply \eqref{e:CH-for-v} to the one dimensional integral appearing in \eqref{e:transformed-integral}, we obtain the upper bound 
$\vert I \vert  \, \sup_{\eta_j\in\RR}|k(\eta_\L)|$. Assume for the moment that $f$ is continuously 
differentiable. Then $\sup\limits_{\eta_j\in\RR}|k(\eta_\L)|\le \int_\RR \left  | \frac{\partial k(\eta_\L)}{\partial \eta_j}\right| d\eta_j$.
Since
\[
 \frac{\partial k(\eta_\L)}{\partial \eta_j} 
= \sum_{k\in\tL_{l+R}}\prod_{p\in\tL_{l+R}} f((B_\L\eta_\L)_l) \frac{\partial }{\partial \eta_j} f((B_\L\eta_\L)_k)
= \sum_{k\in\tL_{l+R}} B_\L(k,j) f'((\omega_\L)_k)\prod_{p\in\tL_{l+R}} f((\omega_\L)_l) 
\]
we can pass back to the $\omega$-variables and establish the bound 
\[
\int_{\RR^{L-1}} d\eta_\L ^{\perp j}  \sup_{\eta_j\in\RR}|k(\eta_\L)|
\le |\det A_\L| \,   \sum_{k\in\tL_{l+R}} |B_\L(k,j)| \, \|f'\|_{L^1}.
\]
Thus \eqref{e:transformed-integral} is bounded by  $ \vert I \vert  \, \|B_\L\|_1 \, \|f\|_{\Var} $.
To extend this estimate to general densities $f$ of bounded variation, 
let $\{f_k\}_k$ be an approximation sequence of smooth, nonnegative, compactly supported functions 
such that $\|f\|_1 =1$ for all $k\in\NN$, $\lim_{k\to\infty} \|f_k\|_{\Var}	 = \|f\|_{\Var}$ and
$\lim_{k\to\infty} \|f_k - f\|_1  = 0$. Then we have
\begin{multline*}
\int_{\RR^L} d\omega_\Lambda\; F(\omega_\Lambda)\ \langle \phi , \chi_j P_\omega^l(I) \chi_j  \phi \rangle
=
\int\limits_{\RR^L} d\omega_\Lambda\; \prod_{k\in {\Lambda}} f_k(\omega_k)\
\langle \phi , \chi_j P_\omega^l(I) \chi_j  \phi \rangle
\\
 + \int_{\RR^L} d\omega_\Lambda\; \left[
\prod_{k\in {\Lambda}} f(\omega_k)-\prod_{k\in {\Lambda}} f_k(\omega_k)\right]\ \langle \phi , \chi_j P_\omega^l(I) \chi_j  \phi \rangle.
\end{multline*}
As we have shown, the first integral on the right is bounded by $ |I| \, \|B\|_1\, \| f_k \|_{\Var}$.
This expression tends to $|I| \, \|B\|_1 \| f \|_{\Var}$ for $k\to\infty$.
A telescoping argument shows that the norm of the second integral is bounded by $ L \, \|f-f_k\|_1$,
which tends to zero as $k\to\infty$.
Thus the proposition is proven.
 \end{proof}

The following lemma will be used in the proof of Proposition \ref{p:main2}.

\begin{lem} \label{l:ALambda}
Let $N\in \NN$ and $\pi=\pi_N$. Then, for all $j \in \tL_{N} $ and $ k \in \tL_{N+D} $ 
the equality $\alpha_{\pi(j-k)} = \alpha_{j-k}$ holds. 
\end{lem}
\begin{proof}
We first show that $\alpha_{\pi(j-k)}=0$ implies $\alpha_{j-k}=0$: Set $m=j-k$.
If $\alpha_m \neq 0$ then  $m \in \L_D$, since the support of $\alpha$ is contained in $\L_D$. 
This implies $\pi(m)=m$ and thus $\alpha_{\pi(m)}=\alpha_{m}$.

Now we consider the case $\alpha_{\pi(m)}\neq 0$. In this case $m$ is an element of $\L_D + \Gamma_{N+D}$.
Since $j \in \tL_{N+r} $ and $ k \in \tL_{N+D} $ , the triangle inequality implies
$\|j-k\|_\infty\le 2N +D$.  As
\[
 \Lambda_{2N+D} \cap (\Lambda_D + \Gamma_{N+D})= \Lambda_D
\]
we conclude that $\pi(j-k)=j-k$.
\end{proof}

\begin{proof}[Proof of Proposition \ref{p:main2}]
An application of Lemma \ref{l:ALambda} with the choice $N=l+r$ 
(recall that $R \le r+D$) shows that the matrix $A_\Lambda$ with coefficients  
$A_{\Lambda}({j,k}) = \alpha_{\pi(j-k)}$ for  $j,k\in\tL_{l+R}$
satisfies condition \eqref{e:Acondition}. Now we have to identify the inverse of $A_\L$.
Define $B_\Lambda \colon \ell^1(\tL_{l+R})\to \ell^1(\tL_{l+R})$ by
$ \displaystyle  B_\Lambda(j,k) := \sum_{p \in \Gamma_{l+R}(j)} B(p,k)$
where $B(p,k)$ are the coefficients of the inverse $B$ of $A\colon \ell^1(\ZZ^d)\to \ell^1(\ZZ^d)$.
Recall that the projection $\pi \colon \ZZ^d \to \Lambda_{l+R}$ is defined by 
$\pi (m)=  \Lambda_{l+R} \,\cap\, \Gamma_{l+R} (m)$.
Let us calculate the product $ B_\Lambda A_\Lambda$. 
For any $i,j \in Q_{l+R}$ we have by the very definition of $A_\L$ and $B_\L$:
\begin{align*}
 \label{e:inverse-calculation}
(B_\Lambda A_\Lambda) (i,j) 
&= \sum_{m\in \tL_{l+R}} B_\L(m,j) \alpha_{\pi(i-m)} 
= \sum_{m\in \tL_{l+R}} \sum_{n\in\Gamma} B(m+n,j)  \sum_{p\in\Gamma} \alpha_{i-m+p} 
\\&= \sum_{m\in \tL_{l+R}} \sum_{n\in\Gamma} B(m+n,j) \sum_{p\in\Gamma} \alpha_{i-m+p-n} .
\intertext{In the last line we used that $n+\Gamma= \Gamma$ for $n\in\Gamma$. 
By the definition of the Laurent matrix $A$, the last expression equals}
&= \sum_{p\in\Gamma} \sum_{k \in \ZZ^d} B(k,j) A(i+p,k)
= \sum_{p\in\Gamma} \delta_{i+p,j} 
\end{align*}
since $B$ is the inverse of $A$. 
Note that for $i,j, \in \Lambda_{l+R}$ we have $\delta_{i+p,j} =0$ for all  $p \in \Gamma \setminus \{0\}$.
It follows that
\[
 \sum_{p\in\Gamma} \delta_{i+p,j} =\delta_{i,j}.
 \]
Thus we have checked that $B_\L$ is the inverse of $A_\L$.

The last step in the proof is to establish that $\|B_\L\|_1 \le \|B\|_1$. Indeed, for all $\L$ we have
\[
\|B_\L\|_1 \le 
 \sum_{j\in\L_{l+R}}  \, | B_\L(j,k)| 
\le  \sum_{j\in\L_{l+R}}  \, | \sum_{p\in \Gamma} B_\L(j+p,k)| 
\le  \sum_{m\in\ZZ^d}  \, |B_\L(j+p,k)| \le  \|B\|_1.
\]
 \end{proof}

\begin{proof}[Completion of the proof of Theorem \ref{t:main}]
We collect the estimates established so far:
\begin{align*}
 \EE &\left \{ \Tr P_\omega^l(I) \right \} 
 \le C(E_2,V)\, \sum_{j\in \tL_l} \left \|  \EE \{\chi_j P_\omega^l(I) \chi_j  \} \right \|
 &\qquad \text{ by Proposition  \ref{p:CombesHislop} (a) }
 \\&
 \le C(E_2,V)\, \sum_{j\in \tL_l} \sup_{\|\phi \|=1} \EE \left \{ \langle \phi , \chi_j P_\omega^l(I) \chi_j \phi \rangle  \right\}
 \\&
 \le C(E_2,V)\, |\tL_l| \,   \|f\|_{\Var} \|B_\Lambda\|_1 \, |I|
 &\qquad \text{ by Proposition  \ref{p:main1}  }
 \\&
 \le C(E_2,V)\, (2l+1)^d \,  \|f\|_{\Var} \|B\|_1 \, |I|
 &\qquad \text{ by Proposition  \ref{p:main2}   }
\end{align*}
 
\end{proof}

\begin{rem}
Daniel Lenz pointed out to us that the calculation in the proof of Proposition \ref{p:main2}
can be understood from a more abstract point of view.  If $G$ is a locally compact abelian group and $H$ a closed subgroup 
(with appropriate invariant measures)
then the projection map induces a homomorphism between the convolution algebras $\ell^1(G)$ and $\ell^1(H)$. 
\end{rem}

\def\cprime{$'$}\def\polhk#1{\setbox0=\hbox{#1}{\ooalign{\hidewidth
  \lower1.5ex\hbox{`}\hidewidth\crcr\unhbox0}}}


\begin{thebibliography}{CHKN02}

\bibitem[Arv02]{Arveson-02}
W.~Arveson.
\newblock {\em A short course on spectral theory}, volume 209 of {\em Graduate
  Texts in Mathematics}.
\newblock Springer-Verlag, New York, 2002.

\bibitem[BS99]{BoettcherS-99}
A.~B{\"o}ttcher and B.~Silbermann.
\newblock {\em Introduction to large truncated {Toeplitz} matrices}.
\newblock Springer, 1999.

\bibitem[CH94]{CombesH-94b}
J.-M. Combes and P.D. Hislop.
\newblock Localization for some continuous, random {Hamiltionians} in
  d-dimensions.
\newblock {\em J. Funct. Anal.}, 124:149--180, 1994.

\bibitem[CHK07]{CombesHK-07}
J.-M. Combes, P.~Hislop, and F.~Klopp.
\newblock An optimal {W}egner estimate and its application to the global
  continuity of the integrated density of states for random {S}chr\"odinger
  operators.
\newblock {\em Duke Math. J.}, 140(3):469--498, 2007.

\bibitem[CHKN02]{CombesHKN-02}
J.-M. Combes, P.~D. Hislop, F.~Klopp, and Shu Nakamura.
\newblock The {W}egner estimate and the integrated density of states for some
  random operators.
\newblock {\em Proc. Indian Acad. Sci. Math. Sci.}, 112(1):31--53, 2002.
\newblock www.ias.ac.in/mathsci/.

\bibitem[GK04]{GerminetK-04}
F.~Germinet and A.~Klein.
\newblock A characterization of the {A}nderson metal-insulator transport
  transition.
\newblock {\em Duke Math. J.}, 124(2):309--350, 2004.
\newblock http://www.ma.utexas.edu/mp\_arc/c/01/01-486.pdf.

\bibitem[HK02]{HislopK-02}
P.~D. Hislop and F.~Klopp.
\newblock The integrated density of states for some random operators with
  nonsign definite potentials.
\newblock {\em J. Funct. Anal.}, 195(1):12--47, 2002.
\newblock www.ma.utexas.edu/mp\_arc, preprint no. 01-139 (2001).

\bibitem[Kat04]{Katznelson-04}
Y~Katznelson.
\newblock {\em An introduction to harmonic analysis}.
\newblock Cambridge Mathematical Library. Cambridge University Press,
  Cambridge, third edition, 2004.

\bibitem[Klo95]{Klopp-95a}
F.~Klopp.
\newblock Localization for some continuous random {Schr\"odinger} operators.
\newblock {\em Commun. Math. Phys.}, 167:553--569, 1995.

\bibitem[KM07]{KirschM-07}
W.~Kirsch and B.~Metzger.
\newblock The integrated density of states for random {Schr\"odinger}
  operators.
\newblock In {\em Spectral Theory and Mathematical Physics}, volume~76 of {\em
  Proceedings of Symposia in Pure Mathematics}, pages 649--698. AMS, 2007.

\bibitem[KS80]{KozakS-80b}
A.~V. Kozak and I.~B. Simonenko.
\newblock Projection methods for solving multidimensional discrete convolution
  equations.
\newblock {\em Sibirsk. Mat. Zh.}, 21(2):119--127, 237, 1980.

\bibitem[KS87]{KotaniS-87}
S.~Kotani and B.~Simon.
\newblock Localization in general one-dimensional random systems. {I}{I}.
  {C}ontinuum {S}chr\"odinger operators.
\newblock {\em Comm. Math. Phys.}, 112(1):103--119, 1987.

\bibitem[KV06]{KostrykinV-06}
V.~Kostrykin and I.~Veseli\'c.
\newblock On the {Lipschitz} continuity of the integrated density of states for
  sign-indefinite potentials.
\newblock {\em Math. Z.}, 252(2):367--392, 2006.
\newblock http://arXiv.org/math-ph/0408013.

\bibitem[Sto01]{Stollmann-01}
P.~Stollmann.
\newblock {\em Caught by disorder: Bound States in Random Media}, volume~20 of
  {\em Progress in Mathematical Physics}.
\newblock Birkh\"auser, 2001.

\bibitem[Ves02]{Veselic-02a}
I.~Veseli\'c.
\newblock {W}egner estimate and the density of states of some indefinite alloy
  type {S}chr\"odinger operators.
\newblock {\em Lett. Math. Phys.}, 59(3):199--214, 2002.
\newblock http://www.ma.utexas.edu/mp\_arc/c/00/00-373.ps.gz.

\bibitem[Ves07]{Veselic-07b}
I.~Veseli\'c.
\newblock {\em {\it Existence and regularity properties of the integrated
  density of states of random {Schr\"odinger} Operators}}, volume Vol. 1917 of
  {\em Lecture Notes in Mathematics}.
\newblock Springer-Verlag, 2007.

\bibitem[Weg81]{Wegner-81}
F.~Wegner.
\newblock Bounds on the {DOS} in disordered systems.
\newblock {\em Z. Phys. B}, 44:9--15, 1981.

\end{thebibliography}
\end{document}